\documentclass{amsart} 
\usepackage[all]{xy}
\usepackage{amsmath,amsfonts,epsfig, psfrag, graphics}

\begin{document} \bibliographystyle{plain}

\newtheorem{thm}{Theorem}[section] \newtheorem{lem}[thm]{Lemma} \newtheorem{cor}[thm]{Corollary}
\newtheorem{mainlem}[thm]{Main Lemma} \newtheorem{prop}[thm]{Proposition} \newtheorem{conj}[thm]{Conjecture}

\theoremstyle{definition} \newtheorem{defn}{Definition}[section]

\newcommand{\U}{\ensuremath{\widetilde}} \newenvironment{pfsketch}{{\it Sketch of Proof:}\quad}{\square \vskip
12pt}

\newtheorem{rmk}[thm]{Remark}

\newcommand{\Hn}{\ensuremath{\mathbb{H}^n}} \newcommand{\h}{\ensuremath{{\text{hyp}}}}
\newcommand{\g}{\ensuremath{{\text{geod}}}} \newcommand{\Hess}{\ensuremath{{\text{Hess} \ }}}
\newcommand{\tr}{\ensuremath{{\text{tr} \, }}} \newcommand{\deri}{\ensuremath{\frac{\partial}{\partial x_i}}}
\newcommand{\derj}{\ensuremath{\frac{\partial}{\partial x_j}}} \newcommand{\V}{\ensuremath{\text{Vol} }}

\def\square{\hfill${\vcenter{\vbox{\hrule height.4pt \hbox{\vrule width.4pt height7pt \kern7pt \vrule width.4pt}
\hrule height.4pt}}}$}

\newenvironment{pf}{{\it Proof:}\quad}{\square \vskip 12pt}

\title{The minimal entropy conjecture for nonuniform rank one lattices}
\author{Peter A. Storm} \date{January 6th, 2006}

\begin{abstract} The Besson-Courtois-Gallot theorem is proven for noncompact finite volume Riemannian manifolds.  In particular, no bounded geometry assumptions are made.  This proves the minimal entropy conjecture for nonuniform rank one lattices.
\end{abstract}

\thanks{This research was partially supported by an NSF Postdoctoral Fellowship.} \maketitle

\section{Introduction}\label{introduction} This paper proves the Besson-Courtois-Gallot theorem  \cite{BCGlong} for finite volume spaces. \begin{thm} \label{main thm} For $n \ge 3$, let $(M,b)$ be a finite volume orientable
Riemannian $n$-manifold.  Let $(Z,g_0)$ be an $n$-dimensional orientable finite volume negatively curved locally
symmetric space.  If $f: M \longrightarrow Z$ is a proper map of nonzero degree then $$h(b)^n \, \text{Vol}(M,b)
\ge | \text{deg}(f) | \, h(g_0)^n \, \text{Vol}(Z,g_0),$$ with equality if and only if $f$ is proper homotopic to
a homothetic covering map. \end{thm}

{\noindent}(See Theorems \ref{inequality thm} and \ref{equality prop}.)  Here $h(b)$ (resp. $h(g_0)$) denotes the volume growth entropy of the Riemannian universal cover of $(M,b)$ (resp. $(Z,g_0)$).  Theorem \ref{main thm} places no bounded geometry assumptions on $(M,b)$.  Results similar to Theorem \ref{main thm} were previously obtained by
Boland-Connell-Souto \cite{BCS, BCS1} under additional hypotheses.

\vskip 6pt The author thanks Richard Canary, Chris Connell, Gilles Courtois, Juha Heinonen, Juan Souto, and Ralf Spatzier for many helpful conversations, and for  listening to various preliminary versions of Theorem \ref{main thm}.  The work of this paper was inspired by \cite{BCGlong} and \cite{BCS}.  

\section{Sketch of the Proof} \label{sketch}

In this section, we assume the reader is familiar with the seminal Besson-Courtois-Gallot paper \cite{BCGlong}.
The proof of Theorem \ref{main thm} differs from the original proof of the Besson-Courtois-Gallot theorem
\cite[pg.734]{BCGlong} in two places.  At these places, the original proof fails in the noncompact setting, and a
new argument is required to arrive at the same conclusion.  Rather than re-sketch arguments from \cite{BCGlong},
we skip straight to the ideas of this paper not found elsewhere.  To streamline the exposition, several extra
assumptions will be made in this section.

Let $(Z,g_0)$ be a noncompact oriented finite volume $n$-dimensional negatively curved locally symmetric space, for $n \ge 3$.  Let $(Z,b)$ be
another finite volume Riemannian metric on $Z$.  As in \cite{BCGlong}, for $c> h(b)$ there exist uniformly Lipschitz $\mathcal{C}^1$-smooth $\pi_1$-equivariant maps 

\begin{eqnarray*} 
\Phi_0  :  (\U{Z}, g_0) & \longrightarrow & S^\infty_+ \subset L^2 (\partial (\U{Z}, g_0))  \\
\Phi_c  : (\U{Z},b) & \longrightarrow &  S^\infty_+ \subset L^2 (\partial (\U{Z}, g_0)), 
\end{eqnarray*} 

{\noindent}where $S^\infty_+$ is the subset of strictly positive functions of unit norm, $\Phi_0$ is the square root of the Poisson kernel, and $\Phi_c$ is a synthetically defined analogue of $\Phi_0$ for a general Riemannian manifold.  If $\text{bar} : L^2_+ (\partial \U{Z}) \longrightarrow (\U{Z}, g_0)$ is the barycenter map, and $\Omega$ is the
closed $n$-form $\text{bar}^* d\text{vol}_{g_0}$, then $\Phi_0$ calibrates the restricted closed $n$-form
$\Omega|_{S^\infty_+}$.  Let $H: \U{Z} \times [0,1] \longrightarrow L^2$ be the straight line homotopy from
$\Phi_0$ to $\Phi_c$.  As everything is equivariant, the closed $n$-forms $\Phi_0^* \Omega, \Phi_c^* \Omega$, and
$H^* \Omega$ are $\pi_1$-invariant and descend to forms on $Z, Z$ and $Z \times [0,1]$ respectively. \\

{\noindent}\textbf{Part 1 }(Sections \ref{homotopy section}-\ref{BCG inequality section})\textbf{:}  \emph{Prove
that $h(b)^n \, \text{Vol}(Z,b) \ge h(g_0)^n \, \text{Vol}(Z,g_0)$.}

In the noncompact setting, the proof of this inequality in \cite{BCGlong} fails at one key point.  Namely, we
would like to claim the equation \begin{equation}\label{sketch 1} \int_Z \Phi_c^* \Omega - \int_Z \Phi_0^* \Omega
\stackrel{?}{=} \int_{Z \times [0,1]} d H^* \Omega = \int_{Z \times [0,1]} 0 =  0 \end{equation} holds by Stokes' Theorem.  But $Z$ is
noncompact, so Stokes' Theorem may not be applied directly.  Instead we must pick a compact exhaustion $Z_i
\subset Z_{i+1} \subset \ldots$ of $Z$ by codimension zero submanifolds, and apply Stokes' theorem to each $Z_i$.
This yields \begin{equation}\label{sketch 2} \int_{Z_i} \Phi_c^* \Omega - \int_{Z_i} \Phi_0^* \Omega \pm
\int_{\partial Z_i \times [0,1]} H^* \Omega = 0. \end{equation} Equation (\ref{sketch 1}) would follow by proving
the error term in equation (\ref{sketch 2}), i.e. the third summand, goes to zero as $i \rightarrow \infty$.  This
estimate involves:\\ $\bullet \ \ $ Finding a good compact exhaustion $\{ Z_i \}_i$ such that the
$(n-1)$-dimensional volume of $(\partial Z_i,b)$ goes to zero.  This is easy using a global smooth proper $1$-Lipschitz function on $(Z,b)$, which always exists \cite{GW}.  \\ $\bullet \ \ $ Finding a map $f: Z
\longrightarrow Z$ properly homotopic to the identity such that $f |_{\partial Z_i}$ is $1$-Lipschitz for all $i$.  This is possible because the ends of the locally symmetric space $(Z,g_0)$ are negatively curved cusps, which contract exponentially toward infinity.  The map $f$ can be pushed out the cusps of $(Z,g_0)$ until it is $1$-Lipschitz on the slices $(\partial Z_i ,b )$.  Homotopies which push maps further out the ends are proper homotopies, so this deformation of $f$ is admissible, despite the fact that it will not necessarily have bounded tracks in $(Z,g_0)$.
\\ $\bullet \ \ $ Controlling the stretching of the straight line homotopy from $\Phi_c$ to $\Phi_0 \circ f$.  This follows from elementary estimates.
\\ After proving the error term in equation (\ref{sketch 2}) goes to zero, equation (\ref{sketch 1}) follows.  The
rest of Part 1 proceeds as in \cite{BCGlong}. \\

{\noindent}\textbf{Part 2 }(Section \ref{equality case})\textbf{:} \emph{If $h(b)^n \, \text{Vol}(Z,b) = h(g_0)^n
\, \text{Vol}(Z,g_0)$ then prove that the identity map is properly homotopic to an isometric homeomorphism.}

Assume \begin{equation}\label{sketch 3} h(b)^n \, \text{Vol}(Z,b)= h(g_0)^n \, \text{Vol}(Z,g_0). \end{equation}
After possibly rescaling the metrics, we may also assume that $h(b) = h(g_0) = h$.  For $c> h$, define the
natural map $(\U{Z},b) \longrightarrow (\U{Z},g_0)$ as the composition of $\Phi_c : (\U{Z},b) \longrightarrow
L^2_+$ and the barycenter map $\text{bar}: L^2_+ \longrightarrow (\U{Z}, g_0) \cong \mathbb{H}^n$.  These maps
descend to a family of maps $F_c : (Z,b) \longrightarrow (Z,g_0)$, which are each homotopic to the identity map $Z \longrightarrow Z$.  Using equation (\ref{sketch 3}) we show that
$\text{Jac} \, F_c \le (c/h)^n$ and $\text{Jac} \, F_c$ converges to the constant function $1$ in $L^1
(Z,d\text{vol}_b)$.

For a generic point $p \in (Z,b)$, we can force some sequence of maps $F_{c_i}$ to converge uniformly on a small
neighborhood of $p$.  This is done by post-composing each $F_c$ by an appropriate orientation preserving isometry
$\gamma_c \in \text{Isom}^+ (\U{Z},g_0)$.  Using the arguments of \cite{BCGlong}, we may conclude that the maps
$\gamma_{c_i} \circ F_{c_i}$ converge uniformly near $p$ to a $1$-Lipschitz map $F^p$.  A short topological argument shows
that these locally defined maps can be pieced together to form a global $1$-Lipschitz map $F: (Z,b)
\longrightarrow (Z,g_0)$.

The remainder of the argument is devoted to proving that $\text{Jac} \, F = 1$ a.e..  If this is not true, then
$F$ strictly shrinks volume by some definite amount on a fixed compact submanifold of $(Z,b)$.  After choosing a submanifold $\{\partial Z_i \}$ which is very far from this compact set, we modify the map $F$ near $\partial Z_i$ to make it proper, while adding very little volume to its image in $(Z,g_0)$.  This modification involves coning off the map $F$ near the slice $\partial Z_i$, forcing it to leave compact sets at a definite rate.  This coning operation adds very little volume to the image because the target manifold $(Z,g_0)$ has pinched negative
curvature.  Therefore even after coning the modified map still globally shrinks volume by some definite amount.  A short topological argument shows that the modified map is proper, and so we can apply the change of variables formula to it.  By the hypothesis that $(Z,b)$ and $(Z,g_0)$ have equal volume, the change of variables formula implies that the modified map \emph{cannot shrink volume globally}.  This yields a
contradiction, showing that $\text{Jac} \, F = 1$ a.e..

We may then apply the theory of quasiregular maps (see Section \ref{quasiregular maps section}) to conclude that
$F$ is in fact a local isometry.  This implies that $(Z,b)$ is a finite volume negatively curved locally
symmetric space.  With this information the rest of the proof follows easily.

\section{Preliminaries}

All Riemannian manifolds in this paper are assumed to be oriented and smooth with smooth Riemannian metrics.

Once and for all, let us fix the notation of Theorem \ref{main thm}:  For $n \ge 3$, let $(M,b)$ be a connected finite volume Riemannian $n$-manifold.  Let $(Z,g_0)$ be an $n$-dimensional finite volume negatively curved locally symmetric space.  Let $f: M \longrightarrow Z$ be a proper map of nonzero degree.  Since only the proper homotopy class of $f: M \longrightarrow Z$ is important, we may assume without a loss of generality that $f$ is smooth.  We will denote the Riemannian universal cover of $(M,b)$ (resp. $(Z,g_0)$) by $(\U{M},b)$ (resp. $(\U{Z},g_0)$).  Fix a basepoint $o \in \U{Z}$.  Let $f: \U{M} \longrightarrow \U{Z}$ also denote the lifted map between universal covers.  (To avoid confusion, the domain and range of the map will always be specified.)

\subsection{Volume growth entropy} Let $(N,g)$ be a connected Riemannian manifold.  Let $({\U{N}},g)$ be the universal cover
of $N$ equipped with the lifted metric.  The \textit{volume growth entropy} of $({\U{N}},g)$ is the number $$h(g)
:=  \limsup_{R \rightarrow \infty}
   \frac{1}{R} \log \text{Vol}_g ( B_{{\U{N}}} (x,R)),$$
where $x \in {\U{N}}$, and $B_{{\U{N}}} (x,R) \subset \widetilde{N}$. \vskip 3pt {\noindent}(The volume growth
entropy is independent of the choice of $x \in {\U{N}}$.)

\subsection{Visual measure} For the negatively curved rank one symmetric space $(\U{Z},g_0)$, let $P: S_p \U{Z} \longrightarrow \partial \U{Z}$ be the standard radial homeomorphism between the unit tangent sphere and the boundary
at infinity.  Normalize Lebesgue measure on $S_p \U{Z}$ to have mass $1$.  Define $\mu_p$ to be the push-forward by
$P$ of the normalized Lebesgue measure on $S_p \U{Z}$.  This probability measure is known as \textit{visual measure
at} $p$.

\subsection{The barycenter map} \label{the barycenter map} Let $\{B_\theta \}_{\theta \in \partial Z}$ be the Busemann
functions of $(\U{Z},g_0)$, normalized such that $B_\theta (o) =0 $ for all $\theta \in \partial \U{Z}$.  Consider
the Hilbert space $L^2 (\partial \U{Z})$ of square integrable functions on $\partial \U{Z}$ with respect to the
visual measure $\mu_o$.  Define an Isom$(\U{Z},g_0)$-action on $L^2 (\partial \U{Z})$ by $$(\gamma. \phi)(\theta)
:= \phi(\gamma^{-1}.\theta) \cdot \sqrt{\exp{(-h(g_0) B_\theta (\gamma.o))}}.$$ Then Isom$(\U{Z},g_0)$ acts by
isometries on $L^2 (\partial \U{Z})$.  Let $L^2_+$ denote the strictly positive functions in $L^2 (\partial
\U{Z})$.  Note that Isom$(\U{Z},g_0)$ acts by isometries on $L^2_+$.

Following \cite[Sec.5]{BCGlong}, define the barycenter map $\text{bar} : L^2_+   \longrightarrow   (\U{Z},g_0) $
by the implicit equation $$\int_{\partial \U{Z}} \langle \nabla B_\theta \, , \, v \rangle_{\text{bar}(\phi)}\
\phi(\theta)^2 d \mu_o (\theta) =0$$ for all $v \in T_{\text{bar}(\phi)} \U{Z}$.  The barycenter map is well
defined, $\mathcal{C}^1$-smooth, and $\text{bar} (c \phi) = \text{bar} (\phi)$ for all $\phi \in L^2_+$ and $c>0$
\cite[Sec.5]{BCGlong}.

\subsection{Quasiregular maps} \label{quasiregular maps section} The following fact will be used in the proof of
Proposition \ref{equality prop}. \begin{thm} \label{qr maps thm} Let $U\subset \mathbb{R}^n$ be open.  Let
$(U,g_U), (\mathbb{R}^n, g)$ be Riemannian manifolds.  If $f: (U,g_U) \longrightarrow (\mathbb{R}^n, g)$ is
locally $1$-Lipschitz and $\text{Jac} \, f =1$ almost everywhere in $U$, then $f$ a local homeomorphism and a
local isometry. \end{thm}
{\noindent}By the 1970s, this fact was well known to many Siberian and Finnish analysts. 
In particular, Theorem \ref{qr maps thm} follows
from the more general \cite[Thm.10.5]{Re}, proved independently by Gol'dshtein and Martio-Rickman-V{\"a}is{\"a}l{\"a} (both building on work of Reshetnyak \cite{Re}).  (\cite[Thm.10.5]{Re} is stated only for Euclidean space.  As usual,
the Riemannian case follows by considering small neighborhoods nearly isometric to Euclidean neighborhoods.)  A
self-contained proof of Theorem \ref{qr maps thm} can be obtained from the arguments of \cite[App.C]{BCGlong}.  The statements of
\cite[App.C]{BCGlong} are made for closed manifolds, but by using local topological degree the arguments can be
applied locally to prove Theorem \ref{qr maps thm}.

\section{The map into Hilbert space} \label{map into Hilbert space section} 

This section is a restatement of facts  from \cite{BCGlong}.  Following \cite[Sec.2.3,Sec.8]{BCGlong}, in this section we define a family of $\pi_1$-equivariant $\mathcal{C}^1$
maps $$\Phi^b_c : ({\U{M}}, b) \longrightarrow S^\infty_+ \subset L^2 (\partial {\U{Z}}) \text{  for } c> h(b):=
h(({\U{Z}},b)).$$ As before, $L^2 (\partial {\U{Z}})$ is the Hilbert space of square integrable 
functions on $\partial {\U{Z}}$ with respect to the visual measure at the basepoint $o \in {\U{Z}}$.  $S^\infty_+
\subset L^2 (\partial {\U{Z}})$ is the subset of strictly positive functions of unit norm.

For $c > h(b)$ define \begin{eqnarray*} \Psi^b_c : ({\U{M}},b) \times \partial {\U{Z}} & \longrightarrow &
\mathbb{R} \\ (p, \theta) &\longmapsto & \left( \int_{\U{M}} e^{-c d_b (p,y)} e^{-h(g_0) B_\theta (f(y))}
              d \text{vol}_b (y) \right)^{1/2},
\end{eqnarray*} and $$\Phi^b_c (p,\theta) := \frac{ \Psi^b_c (p,\theta)}{ \left( \int_{\partial {\U{Z}}} (\Psi^b_c
(p,\theta))^2 d \theta \right)^{1/2} }.$$ $\Phi^b_c$ defines a map into function space by fixing its first
coordinate: \begin{eqnarray*} \Phi^b_c : ({\U{M}},b) & \longrightarrow & \{ \text{fcts. on } \partial {\U{Z}} \}
\\
        p        & \longmapsto &  \Phi^b_c (p, \cdot).
\end{eqnarray*}

\begin{lem} \cite[Lem.2.4]{BCGlong} \label{hilbert map lemma} The map $p \longmapsto \Phi^b_c (p, \cdot)$ is a
$\pi_1$-equivariant $\mathcal{C}^1$ map from ${\U{M}}$ into $S^\infty_+ \subset L^2 (\partial {\U{Z}})$.
Moreover, for $u \in T_p {\U{M}}$, $$\int_{\partial {\U{Z}}} | u. \Phi^b_c (x, \theta) |^2 \ d\theta \le
\frac{c^2}{4} \ b(u,u)_p$$ (see \cite[pg.430]{BCGold}). \end{lem}

Using $\Phi^b_c$, pull back the metric on $S^\infty_+$ to ${\U{M}}$.  This defines the positive semi-definite
$(0,2)$-tensor $$g_{\Phi^b_c}(u,v) := \langle d \Phi^b_c \, (u), d \Phi^b_c \, (v) \rangle_{L^2 (\partial
{\U{Z}})}.$$ From the inequality of Lemma \ref{hilbert map lemma}, we obtain \cite[pg.744]{BCGlong}
\begin{equation} \label{star} \| d\Phi^b_c (u) \|^2 = g_{\Phi^b_c} (u,u) \le  \frac{c^2}{4} \, b (u,u).
\end{equation} The $\pi_1$-equivariance of $\Phi^b_c$ implies the tensor $g_{\Phi^b_c}$ is $\pi_1$-invariant.  It
thus descends to a tensor $g_{\Phi^b_c}$ on $M$.

\begin{prop} \label{volume inequality} \cite[Prop.3.4]{BCGlong} $$| \text{Jac} \, \Phi^b_c | \le
\left(\frac{c^2}{4n} \right)^{n/2}. \quad \text{In particular,} \quad \text{Vol} (M, g_{\Phi^b_c}) \le \left(
\frac{c^2}{4n} \right)^{n/2}  \text{Vol} (M, b).$$ \end{prop} {\noindent}(Notice inequality (\ref{star}) does not
immediately imply Proposition \ref{volume inequality}.  Better estimates are required.)

\vskip 6pt For the locally symmetric space $(Z,g_0)$, stronger statements are true.  Define \begin{eqnarray*} \Phi_0 : ({\U{Z}}, g_0) &
\longrightarrow &  S^\infty_+ (\partial {\U{Z}}) \\
        p   & \longmapsto &  \exp{( - \frac{h(g_0)}{2} B_\theta (p))} .
\end{eqnarray*} Using $\Phi_0$, pullback the metric on $S^\infty_+$ to ${\U{Z}}$.  This defines the tensor
$$g_{\Phi_0} (u,v) := \langle d \Phi_0 \, (u), d \Phi_0 \, (v) \rangle_{L^2 (\partial {\U{Z}})}.$$ It has been
proven that $$g_{\Phi_0} = \frac{h(g_0)^2}{4n} g_0, \quad \text{  implying  } \quad
    \text{Vol} (Z, g_{\Phi_0}) = \left( \frac{h(g_0)^2}{4n} \right)^{n/2} \text{Vol} (Z, g_0)$$
\cite[pgs.743-744]{BCGlong}.

\section{Calibration theory} \label{calibration theory section}

This section introduces the terminology from calibration theory necessary to properly cite
\cite[Sec.4-5]{BCGlong}.

For the purpose of giving clear definitions, let us for the moment consider a more general situation.  Let $N$ be
an oriented $n$-manifold.  Let $\mathcal{H}$ be a (real) Hilbert manifold (e.g. $S^\infty_+ \subset L^2$).  Assume
there exists an isometric group action $\pi_1 (N) \times \mathcal{H} \longrightarrow \mathcal{H}$.  Any $\pi_1
(N)$-equivariant $\mathcal{C}^1$ map  defines a $\pi_1 (N)$-invariant
symmetric positive semi-definite $(0,2)$-tensor $$g_\Theta (u,v) : = \langle d\Theta \, (u), d \Theta \, (v)
\rangle_{\mathcal{H}}.$$ This tensor descends to $N$, and determines a volume form $d \text{vol}_\Theta$.  Thus we
may define $$\text{Vol}(\Theta) := \int_N d \text{vol}_\Theta.$$

Define the \emph{comass} of an $n$-form $\alpha$ on $\mathcal{H}$ to be $$\text{comass} (\alpha) := \sup |
\alpha_x (v_1, \ldots, v_n) |,$$ where the supremum is taken over all $x \in \mathcal{H}$ and all orthonormal
$n$-tuples $\{ v_i \}_{i=1}^n$ of $\mathcal{H}$.  Then for any $\pi_1 (N)$-invariant $n$-form $\alpha$, we have
the inequality \begin{equation} \label{calibration 1} \text{Vol}(\Theta) \ge \frac{1}{\text{comass}(\alpha)}
\int_N \Theta^* \alpha. \end{equation}
A closed $\pi_1 (N)$-invariant $n$-form $\omega$ \emph{calibrates} an \emph{immersion} $\Theta$ if for all $p \in
N$ and positively oriented bases $\{ v_1, \ldots, v_n \}$ of $T_p N$ we have $$  \omega_{\Theta(p)} (d \Theta
(v_1), \ldots, d \Theta (v_n))  = \text{comass}(\omega) \
            \left( \det{( \langle d\Theta (v_i), d \Theta (v_j) \rangle )_{ij} }\ \right)^{1/2}.$$
The most important property of a calibrating pair $\omega, \Theta$ is
\begin{equation} \label{calibration 2} \text{Vol}(\Theta) = \frac{1}{\text{comass}(\Theta)} \int_N \Theta^*
\omega. \end{equation}

A key tool for proving this paper's main result is \begin{thm} \label{calibration thm} \cite[Prop.5.7]{BCGlong}
Recall the immersion $\Phi_0: {\U{Z}} \longrightarrow S^\infty_+$ of Section \ref{map
into Hilbert space section}, and the $\mathcal{C}^1$ barycenter map $\text{bar}: L^2_+ \longrightarrow
\U{Z}$ of Section \ref{the barycenter map}.  The form $\Omega := \text{bar}^* (d \text{vol}_{g_0})$ is a closed
$\pi_1 (Z)$-invariant $n$-form.  The restricted form $\Omega|_{S^\infty_+}$ calibrates the immersion $\Phi_0$ in
the Hilbert manifold $S^\infty_+$.  Moreover, $$\text{comass} (\Omega|_{S^\infty_+}) = \left( \frac{4n}{h(g_0)^2} \right)^{n/2}.$$
\end{thm}

Recall that $L^2_+ \subset L^2 (\partial \U{Z})$ is the subset of strictly positive functions.  For later reference, we note \begin{lem} \label{bounded comass} If $B(0,.5)$ is the open ball of radius $.5$ about
the origin in $L^2 (\U{Z})$, then the comass of the form $\Omega:= \text{bar}^* (d \text{vol}_{g_0})$ is uniformly
bounded above on the subset $L^2_+ - B(0,.5)$. \end{lem} \begin{pf} Radial projection to $S^\infty_+$ is
$2$-Lipschitz on $L^2_+ - B(0,.5)$.  The barycenter map is invariant under radial projection (see Section \ref{the
barycenter map}).  The lemma now follows from the fact that the comass of $\Omega$ restricted to $S^\infty_+$ is
finite. \end{pf}

\section{Estimates on the straight line homotopy} \label{homotopy section}

Recall that $(Z, g_0)$ is a connected $n$-dimensional negatively curved locally symmetric space, $(M,b)$ is a connected Riemannian
$n$-manifold, and $f: M \longrightarrow Z$ is a proper map of nonzero degree.  Let $f: \U{M} \longrightarrow
\U{Z}$ also denote the lifted map between universal covers.  Recall also that $S^\infty_+ \subset L^2 (\partial
\U{Z})$ is the subset of strictly positive functions of unit norm, and $L^2_+ \subset L^2 (\partial \U{Z})$ is the subset of strictly positive functions.

Let $\Theta, \Upsilon: {\U{M}} \longrightarrow S^\infty_+ \subset L^2 (\partial \U{Z})$ be any two $\mathcal{C}^1$ maps which are equivariant with respect to the $\pi_1 (M)$-action on $S^\infty_+$ induced by the homomorphism $f_* : \pi_1 (M) \longrightarrow \pi_1 (Z)$.  The goal of this section is to prove the straight line homotopy from
$\Theta$ to $\Upsilon$ does not stretch too much.  These estimates will be used in the proof of Theorem
\ref{inequality thm}.

Define the straight line homotopy as \begin{eqnarray*} H : {\U{M}} \times [0,1] \times \partial {\U{Z}} &
\longrightarrow &  \mathbb{R}_+  \\
                (x,t, \theta) & \longmapsto   &
    \left[  (1-t) \, \Theta (x,\theta) + t \, \Upsilon (x,\theta) \right].
\end{eqnarray*} By fixing the first two coordinates, this defines a $\mathcal{C}^1$ map $H : {\U{M}} \times [0,1]
\longrightarrow L^2_+$.  $\pi_1$-equivariance follows from the equivariance of $\Theta, \Upsilon$.

For a slice of the homotopy \begin{eqnarray*} H^t : {\U{M}} & \longrightarrow  &  L^2_+ \\
     x & \longmapsto &  \left\{ \theta \mapsto H(x,t,\theta) \right\} \in L^2_+,
\end{eqnarray*} we now estimate how much $H^t$ stretches ${\U{M}}$ in terms of $\Theta, \Upsilon$.  Let $v \in T_x
{\U{M}}$, $t \in [0,1]$.  Note that $dH^t_x (v)$ is an element of the tangent space $T_{H^t (x)} L^2 (\partial \U{Z}) = L^2
(\partial {\U{Z}})$.  In particular, $dH^t_x (v)$ is a function of $\theta \in \partial {\U{Z}}$.  We compute $$[
dH^t_x (v) ](\theta) =  (1-t) \,
    [ d \Theta_x (v) ] (\theta) \ + \ t \,  [d \Upsilon_x (v) ] (\theta) .$$
Using the inequality $$\left[ (1-t) \delta + t \eta \right]^2 \le (1-t) \delta^2 + t \eta^2 \text{ for all }
  t \in [0,1], \ \delta,\eta \in \mathbb{R},$$
we estimate \begin{eqnarray*} \| dH_x^t (v) \|^2 &= & \int_{\partial {\U{Z}}} [dH_x^t (v)](\theta)^2  \ d\theta
\\ &\le &   \int_{\partial {\U{Z}}} \left\{(1-t) [d \Theta_x (v) ](\theta)^2 + t \, [d\Upsilon_x (v)](\theta)^2
\right\}  \ d \theta  \\
    &\le &    \| d \Theta_x (v) \|^2 + \| d \Upsilon_x (v) \|^2.
\end{eqnarray*}

We now estimate how much $H$ stretches in the time direction.  Let $\frac{\partial}{\partial \tau}$ be the
standard positive basis vector of $T_t [0,1]$.  We compute $$dH_{(x,t)} (\frac{\partial}{\partial \tau}) =
\frac{\partial}{\partial \tau} H (x,\tau) |_{\tau = t}
    =  - \Theta (x) + \Upsilon(x).$$
Suppressing unnecessary notation, we estimate $$ \| \frac{\partial}{\partial \tau} H \|^2  =  \int_{\partial
{\U{Z}}}
    [ \Upsilon - \Theta]^2
    \le  \int_{\partial {\U{Z}}} \left[ \Upsilon^2 + \Theta^2 \right]  = 2.$$
Here we used the fact that $\Upsilon$ and $\Theta$ are maps to positive functions.

\section{The Besson-Courtois-Gallot inequality} \label{BCG inequality section} The goal of this section is to
prove \begin{thm} \label{inequality thm} For $n \ge 3$, let $(M,b)$ be a finite volume oriented Riemannian
$n$-manifold.  Let $(Z,g_0)$ be an $n$-dimensional oriented finite volume negatively curved locally symmetric
space.  If $f: M \longrightarrow Z$ is a proper orientation preserving map of nonzero degree then $$h(b)^n \,
\text{Vol}(M,b) \ge | \text{deg}(f) | \, h(g_0)^n \, \text{Vol}(Z,g_0).$$ \end{thm}

{\noindent}(In the important special case where $M = Z$ and $f: M \longrightarrow Z$ is the identity map, several
technical details do not arise.  Until the argument is clear, the reader is encouraged to read the proof with this
special case in mind.)

For $M$ and $Z$ compact, Theorem \ref{inequality thm} was first proven in \cite{BCGlong}.  If $M$ is not compact
then $Z$ is not compact, because $f$ is proper.  On the other hand, if $Z$ is not compact then $M$ is not compact
because $f$ has nonzero degree.  We may therefore assume without a loss of generality that neither $M$ nor $Z$ is
compact.  Recall ${\U{Z}}$ (resp. $\U{M}$) is the universal cover of $Z$ (resp. $M$).  After possibly reversing the orientation of $M$, we may assume that $f$ has positive degree.  Let $f: \U{M} \longrightarrow \U{Z}$ also denote the lifted map.  By choosing a smooth approximation, we may assume without a loss of generality that $f$ is smooth.

$(Z, g_0)$ is a finite volume negatively curved locally symmetric space.  By the thick-thin decomposition, $Z$ has
a compact submanifold whose complement splits as a smooth manifold into a product region $Y^{n-1}
\times [0,\infty),$ where the manifold $Y$ may not be connected.  The product structure $Y \times [0,
\infty)$ can be chosen to be a standard horospherical foliation.  So we may assume the $(n-1)$-dimensional
$g_0$-volume of the leaf $Y \times \{ r \}$ decreases monotonically to zero as $r \rightarrow \infty$.  Moreover,
we may assume that the shift map \begin{eqnarray*} \sigma_r: Y \times [0, \infty) & \longrightarrow & Y \times [0,
\infty) \\
                         (y, t) & \longmapsto &  (y, t+ r)
\end{eqnarray*} is $K(r)$-Lipschitz with $K(r)$ going to zero monotonically as $r \rightarrow \infty$.

Since Vol$(M, b)$ is finite, we can prove \begin{lem}\label{finding small volume slices} (1)  There exists a
compact exhaustion $\{ M_i \}$ of $M$ such that each $M_i \subset M$ is a smooth submanifold and the
$(n-1)$-dimensional $b$-volume of $\partial M_i$ goes to zero as $i \rightarrow \infty$. \\ (2)  The map $f: M
\longrightarrow Z$ may be altered by a proper homotopy such that the restriction of $f:M \longrightarrow Z$ to
$\cup \partial M_i$ is locally $1$-Lipschitz.
\end{lem}

\begin{pf} There exists a smooth $1$-Lipschitz proper function $\delta: M \longrightarrow [0, \infty)$ (see for
example \cite[Sec.2]{GW}).  This yields the estimate \begin{eqnarray*} \text{Vol}(M,b) \ge \int_M | \nabla \delta
| \, d \text{vol}_b
 & \ge & \int_{\delta^{-1}(\text{reg. val.})} | \nabla \delta | \, d \text{vol}_b  \\
 & = &  \int_{\text{reg. val.}} \left\{ (n-1) \text{ dim'l } b-\text{vol. of } \delta^{-1}(t) \right\} \, dt,
\end{eqnarray*} where ``reg. val.'' indicates the regular values of the function $\delta$.  By Sard's Theorem
there must be a sequence $t_i \rightarrow \infty$ such that $$\left\{ (n-1) \text{ dim'l } b-\text{vol. of }
\delta^{-1}(t_i) \right\} \longrightarrow 0.$$ Define $M_i := \delta^{-1} ([0, t_i])$.  This proves part (1).

Without a loss of generality, we may assume that for all $i$ the image $f(\partial M_i)$ is contained in the product region $Y \times (0, \infty) \subset Z$.  Suppose that for all $i< i_0$, the restricted map $f|_{\partial M_i}$ is $1$-Lipschitz.  We can push
$f(\overline{M - M_{i_0}} )$ out the ends of $(Z,g_0)$ by a family of shift maps $\sigma_r$ until the restricted
map $(\sigma_r \circ f)|_{\partial M_{i_0}}$ is $1$-Lipschitz.  This family of  shift maps can be smoothly
extended to all of $M$ without altering the map $f$ on $M_{i_0 -1} \subset M$, thus yielding a homotopy of $f$. As
this homotopy only pushes further out the ends of $Z$, it is a proper homotopy.  This proves part (2). \end{pf}
{\noindent}(If $M=Z$ then the above proof shows that the identity map is proper isotopic to a diffeomorphism $f:Z
\longrightarrow Z$ such that the $f|_{\partial M_i}$ is $1$-Lipschitz.  In particular, after pulling back the
metric, one may again assume that the identity map satisfies the conclusions of the above lemma.)  Without a loss
of generality, less us assume the proper map $f: M \longrightarrow Z$ satisfies the conclusions of Lemma
\ref{finding small volume slices}.

Recall the two $\pi_1$-equivariant $\mathcal{C}^1$ maps $$\Phi^b_c, \Phi_0 \circ f : {\U{M}} \longrightarrow
S^\infty_+ \subset L^2 (\partial \U{Z}),$$ for $c> h(b)$.  Let $H: {\U{M}} \times [0,1] \longrightarrow L^2_+$ denote the straight line homotopy
(which is $\pi_1$-equivariant) from $\Phi^b_c$ to $\Phi_0 \circ f$.

By Besson-Courtois-Gallot's Theorem \ref{calibration thm}, $\Omega$ is a closed $\pi_1 (Z)$-invariant $n$-form on
$L^2_+$.  When $\Omega$ is restricted to $S^\infty_+$, it calibrates the immersion $\Phi_0$.  $H^* \Omega$ is a
closed $\pi_1 (M)$-invariant $n$-form on ${\U{M}} \times [0,1]$.  It therefore descends to a closed $n$-form on $M
\times [0,1]$.  Moreover, $$H^* \Omega |_{{\U{M}} \times \{0 \}} = \Phi^{b*}_c \Omega \quad \text{ and } \quad
        H^* \Omega |_{{\U{M}} \times \{ 1 \} } =  (\Phi_0 \circ f)^* \Omega.$$
At this stage, we would like to apply Stokes' Theorem to the closed form $H^* \Omega$ on $M \times [0,1]$.  As $M$
is not compact, we are instead forced to apply Stokes' Theorem repeatedly to the compact exhaustion
$M_i \times [0,1] \subset M_{i+1} \times [0,1] \subset \ldots \subset M \times [0,1]$, and then prove the
resulting error terms go to zero.  This will now be carried out.

Define the sequence of submanifolds $L_i := \partial M_i \subset M$.  We may now apply Stokes' Theorem to the
closed form $H^* \Omega$ restricted to $M_i \times [0,1]$.  This yields $$\int_{M_i} {\Phi^b_c}^* \Omega -
\int_{M_i} (\Phi_0 \circ f)^* \Omega = \pm \int_{ L_i \times [0,1]}
            H^* \Omega.$$

\begin{prop} \label{mainprop} $$\lim_{i \rightarrow \infty} \left| \int_{L_i \times [0,1]} H^* \Omega \right| =
0.$$ \end{prop} \begin{pf} Put the product Riemannian metric $b \times (1)$ on $M \times [0,1]$.
$\Phi^b_c$ and $\Phi_0 \circ f$ have their images in positive functions of unit norm in Hilbert space.  By
projecting any pair of points in the image to the plane they span, it follows that the image of the straight line
homotopy $H$ lies in the complement of the ball $B_{L^2} (0,.5)$ of radius $.5$ in $L^2 (\partial \U{Z})$.  By
Lemma \ref{bounded comass}, the comass of $\Omega$ is uniformly bounded on the subset $L^2_+ - B_{L^2} (0,.5)$. By
Lemma \ref{finding small volume slices}, the $(n-1)$-dimensional $b$-volume of the slices $L_i$ goes to zero.
Therefore to prove the proposition it suffices to show that the restricted map $$H|_{L_i \times [0,1]}: (L_i \times [0,1], b \times (1)) \longrightarrow
L^2_+ - B_{L^2}(0,.5)$$ is uniformly Lipschitz for all $i$.  By the estimates of Section
\ref{homotopy section}, at any $(x,t) \in M \times [0,1]$ we have the inequalities $$\| dH^t (v) \|^2 \le \| d\Phi^b_c (v) \|^2 +
\| d (\Phi_0 \circ f) (v) \|^2 \quad \text{and} \quad
    \| \frac{\partial}{\partial \tau} H \|^2  \le 2,$$
where $v \in T_x M$ and $ \frac{\partial}{\partial \tau}$ is a unit basis of $T_t [0,1]$. From inequality
(\ref{star}) of Section \ref{map into Hilbert space section} we have $$ \| d\Phi^b_c (v) \|^2 \le \frac{c^2}{4}
b(v,v) \quad \text{and} \quad
         \| d (\Phi_0 \circ f) (v) \|^2 =  \frac{h(g_0)^2}{4n} g_0 (df(v) \, ,\, df(v)).$$
The restricted map $f: (L_i, b) \longrightarrow (Z, g_0)$ is $1$-Lipschitz by Lemma \ref{finding small volume
slices}.  If $(x,t) \in L_i \times [0,1]$ and $w \in T_x L_i$ then $g_0 (df(w) \, , \, df(w)) \le b(w,w).$ This
implies that $$\| dH^t (w) \|^2 \le \frac{c^2+ h(g_0)^2}{4} \, b(w,w),$$ Therefore $H|_{L_i \times [0,1]}$ is
$\left( \frac{c^2 + h(g_0)^2}{4} + 2 \right)^{1/2}$-Lipschitz, where the ``$+2$" comes from the $\frac{\partial}{\partial \tau}$
direction. \end{pf}

Proposition \ref{mainprop} immediately implies $$\int_M {\Phi^b_c}^* \Omega - \int_M (\Phi_0 \circ f)^* \Omega =
\lim_{i \rightarrow \infty} \left[
            \int_{M_i} {\Phi^b_c}^* \Omega - \int_{M_i} (\Phi_0 \circ f)^* \Omega  \right] = 0 .$$
Recall that by Theorem \ref{calibration thm} the closed $n$-form $\Omega$ restricted to $S^\infty_+$ calibrates
the immersion $\Phi_0$.  Combining this with inequality (\ref{calibration 1}) and equation (\ref{calibration 2}) of Section \ref{calibration theory section} yields \begin{eqnarray*} \text{comass} (\Omega|_{S^\infty_+}) \cdot
\text{Vol}(\Phi^b_c) & \ge & \int_M {\Phi_c^b}^* \Omega =
        \int_M (\Phi_0 \circ f)^* \Omega \\
    &= & \text{deg}(f) \cdot \int_Z \Phi_0^* \Omega =
        \text{deg}(f) \cdot \text{comass} (\Omega|_{S^\infty_+}) \cdot \text{Vol} (\Phi_0).  
\end{eqnarray*} From Section \ref{map into Hilbert space section} we know $$\left( \frac{c^2}{4n} \right)^{n/2}
\text{Vol} (M, b) \ge \text{Vol} (\Phi^b_c) \ \ \text{and} \ \
  \text{Vol} (\Phi_0) = \left( \frac{h(g_0)^2}{4n} \right)^{n/2} \text{Vol} (Z, g_0).$$
This proves that for all $c> h(b)$, $$\left( \frac{c^2}{4n} \right)^{n/2} \text{Vol} (M, b) \ge \text{deg}(f) \
\left( \frac{h(g_0)^2}{4n} \right)^{n/2} \text{Vol} (Z, g_0).$$ Taking $c \searrow h(b)$ yields the desired
inequality, namely $$h(b)^n \, \text{Vol}(M,b) \ge \text{deg}(f) \ h(g_0)^n \, \text{Vol}(Z,g_0).$$ This completes
the proof of Theorem \ref{inequality thm}.

\section{Cones in negatively curved manifolds} \label{cone section}

This section will recall some elementary facts about cones in negatively curved Riemannian manifolds.  These facts will be used in Section \ref{equality case}.

Let $(\U{Z}, g_0)$ be a simply connected $n$-dimensional symmetric space with curvature at most $-1$.  (Everything in this section remains true for $\U{Z}$ any simply connected Riemannian manifold with curvature at most $-1$.)  Let $N$ be a smooth $(n-1)$-dimensional manifold.  Let $g$ be a Riemannian metric on the smooth product $N \times [0, \varepsilon)$.  (The metric $g$ is neither necessarily complete nor a product metric.)  Let $\phi : N \longrightarrow \U{Z}$ be a locally Lipschitz map, where $N$ inherits a Riemannian metric as a submanifold of $(Y \times [0,\varepsilon), g)$.  Let $\theta$ be a point on the boundary at infinity of $\U{Z}$.

With this notation we can define the cone
$${\mathcal{C}} : N \times [0,\varepsilon) \longrightarrow \U{Z}$$
by sending $(x,s) \in N \times [0, \varepsilon)$ to the point at distance $\tan (\frac{s \pi }{2 \varepsilon})$ along the geodesic ray from $\phi(y)$ to $\theta$.

The following elementary lemma is a standard fact whose proof we include for completeness.

\begin{lem} \label{cone lemma}
For any measurable subset $U \subseteq N$ we have the inequality
\begin{equation} \label{cone inequality}
 \int_{U \times [0, \varepsilon)}  | \text{Jac} \, {\mathcal{C}} | \ d \text{vol}_g  \le \frac{1}{n-1} \, 
		\int_U | \text{Jac} \, \phi | \ d \text{vol}_g.
\end{equation}
\end{lem}

{\noindent}Recall that $\phi$ and $\mathcal{C}$ are locally Lipschitz, implying by Rademacher's theorem that they are differentiable almost everywhere.

\begin{pf}
By the change of variables theorem, the value of the integrals in inequality (\ref{cone inequality}) are independent of the coordinate system and the Riemannian metric $g$ on $N \times [0, \varepsilon)$.  So without a loss of generality we are free to make the following changes.  First, we will work on the diffeomorphic manifold $N \times [0, \infty)$, where $\mathcal{C}(x,s)$ is redefined to be the point at distance $s$ along the geodesic from $\phi(x)$ to $\theta$.  Second, we will use $g$ to denote a fixed Riemannian metric $g$ on $N$.  Finally, we equip $N \times [0,\infty)$ with the product metric $g \times (1)$.

At each point where $d \phi$ is injective choose a basis $\{ v_1, v_2, \ldots , v_{n-1} \}$ for $T_x N$ such that the set of vectors $\{ d\phi(v_1), \ldots, d\phi(v_{n-1}) \}$ are orthonormal in $T_{\phi(x)} \U{Z}$.  (The regularity of the resulting vector field is not important.)  Extend each $v_i$ to the vector field $V_i := (v_i, 0)$ on $N \times [0, \infty)$.

For a smooth path $\alpha$ in $N$ such that $\alpha(0) = x$ and $ \alpha' (0) = v_i$, the surface
$$ (t,s) \longmapsto \mathcal{C}(\alpha(t), s)$$
describes a $1$-parameter family of geodesics in $\U{Z}$.  Therefore the vector field $d\mathcal{C}_{(x,s)} (V_i)$ is a Jacobi field along the geodesic
$$ s \longmapsto \mathcal{C}(\alpha(0), s).$$

Let $\pi: T \U{Z} \longrightarrow T\U{Z}$ map a vector in $T_z \U{Z}$ to its component orthogonal to the geodesic from $z$ to $\theta \in \partial \U{Z}$.  As the curvature of $\U{Z}$ is bounded above by $-1$, the Rauch comparison theorem tells us that
$$ \| ( \pi \circ d\mathcal{C} ) (V_i) \|_{\mathcal{C}(x,s)} \le
		e^{-s} \cdot \| ( \pi \circ d\phi) (v_i) \|_{\phi(x)}
			\le e^{-s}.$$
Recall here that $s$ is by definition the distance from $\phi(x) = \mathcal{C}(x,0)$ to $\mathcal{C}(x,s)$.

At points $x \in N$ where $d\phi$ is injective, we obtain the estimate
\begin{eqnarray*}
\lefteqn{\left\| d {\mathcal{C}} \left( V_1 \right) \wedge \ldots \wedge 
					d {\mathcal{C}} \left(V_{n-1} \right) \wedge 
					d {\mathcal{C}} \left( \partial_s \right) \right\|_{\mathcal{C}(x,s)} } \\
 & = & \left\|  (\pi \circ d {\mathcal{C}}) \left( V_1 \right) \wedge \ldots \wedge 
					( \pi \circ d {\mathcal{C}}) \left( V_{n-1} \right) \wedge 
 					d {\mathcal{C}} \left( \partial_s \right) \right\|_{\mathcal{C}(x,s)}  \\
 & \le & e^{- (n-1) s} \cdot 
				\left\| d {\mathcal{C}} \left( \partial_s \right) \right\|_{\mathcal{C}(x,s)} \\
 & = & e^{- (n-1) s} \cdot 
						\left\| d \phi(v_1) \wedge \ldots \wedge d\phi(v_{n-1}) \right\|_{\phi(x)}.
\end{eqnarray*}
Here the first equation comes from the fact that $d\mathcal{C}(\partial_s)$ is orthogonal to the image of $\pi$.  The final equation uses the fact that $d \mathcal{C} (\partial_s)$ is a unit vector and the identity
$\left\| d \phi(v_1) \wedge \ldots \wedge d\phi(v_n) \right\|_{\phi(x)} =1.$

Since $g \times (1)$ is a product metric on $N \times [0, \infty)$, we obtain the equation
$$ \| V_1 \wedge \ldots \wedge V_{n-1} \wedge \partial_s \|_{(x,s)} = \| v_1 \wedge \ldots \wedge v_{n-1} \|_x.$$
This fact yields the following inequality for all $x$ where $d\phi$ is injective.

\begin{eqnarray*}
 | \text{Jac}\, \mathcal{C} (x,s) |  	& = & 
			\frac{ \left\| d \mathcal{C} (V_1) \wedge \ldots \wedge d \mathcal{C}(V_{n-1}) \wedge d \mathcal{C}(\partial_s) \right\|_{\mathcal{C}(x,s)}
					}{ \| V_1 \wedge \ldots \wedge V_{n-1} \wedge \partial_s \|_{(x,s)} }  \\
 & \le &  e^{- (n-1) s}  
				\cdot \frac{ \| d\phi(v_1) \wedge \ldots \wedge d\phi(v_{n-1}) \|_{\phi(x)}}{ \| v_1 \wedge \ldots \wedge v_{n-1} \|_x } \\
 & = & e^{-(n-1) s }  
				\cdot | \text{Jac}\, \phi (x) |.  
\end{eqnarray*}
So finally we obtain the desired estimate
\begin{eqnarray*}
\int_{U \times [0, \infty) } | \text{Jac}\, \mathcal{C} | \, d\text{vol}_{g \times (1)} 
& \le & \int_0^\infty \int_U 
				e^{- (n-1) s } \cdot 				 | \text{Jac}\, \phi (x) | \, d\text{vol}_g (x) \, ds                 \\
& = & \left[ \int_0^\infty 
				e^{- (n-1) s}  				 \, ds \right]
				\cdot \left[ \int_U | \text{Jac}\, \phi | \, d\text{vol}_g \right]   \\
& = & \frac{1}{n-1} \, \int_U | \text{Jac}\, \phi | \, d\text{vol}_g.
\end{eqnarray*}
\end{pf}

We now describe an extension of Lemma \ref{cone lemma} to the case of some negatively curved manifolds with nilpotent fundamental group.  Let $P$ be a group of isometries of $\U{Z}$ which fixes $\theta \in \partial \U{Z}$.  Suppose we are given a locally Lipschitz map $\phi: N \longrightarrow \U{Z}/P$.  Let $\U{N}$ denote the universal cover of $N$, and let 
$\U{\phi} : \U{N} \longrightarrow \U{Z}$ denote the lift of $\phi$.  Notice that the cone
$$\mathcal{C} : \U{N} \times [0, \varepsilon) \longrightarrow \U{Z}$$
is equivariant with respect to the actions of $\pi_1(N)$ on $\U{N} \times [0, \varepsilon)$ and $P$ on $\U{Z}$.  So the ``downstairs'' cone
$$\mathcal{C}: N \times [0, \varepsilon) \longrightarrow \U{Z}/P$$
is well defined.  By applying Lemma \ref{cone lemma} to a fundamental domain for $N$ in $\U{N}$, we obtain the following corollary.

\begin{cor} \label{cone cor}
For any Riemannian metric $g$ on $N \times [0,\varepsilon)$ we have the inequality
$$\int_{N \times [0,\varepsilon)}  | \text{Jac} \, \mathcal{C} | \, d \text{vol}_g 
		\le \frac{1}{n-1} \, \int_N | \text{Jac} \, \phi | \, d \text{vol}_g.$$
\end{cor}

\section{The equality case} \label{equality case} The goal of this section is to prove

\begin{thm} \label{equality prop} For $n \ge 3$, let $(M,b)$ be a finite volume oriented Riemannian $n$-manifold.
Let $(Z,g_0)$ be an $n$-dimensional oriented finite volume negatively curved locally symmetric space.  If $f: M
\longrightarrow Z$ is a proper orientation preserving map of nonzero degree and $$h(b)^n \, \text{Vol}(M,b) =
\text{deg}(f) \, h(g_0)^n \, \text{Vol}(Z,g_0),$$ then $f$ is proper homotopic to a homothetic covering map.
\end{thm}

Since the quantities $h(b)^n \, \text{Vol}(M,b)$ and $h(g_0)^n \, \text{Vol}(Z,g_0)$ are scale invariant, we may assume without a loss of generality that the curvature of $g_0$ is bounded above by $-1$, and
$\text{Vol}(M,b) = \text{deg}(f) \, \text{Vol}(Z,g_0)$.  This implies $h := h(b) = h(g_0).$
For $c> h$, define
the map 
\begin{eqnarray*}
F_c : (M,b) & \longrightarrow & (Z, g_0) \\
        x   & \longmapsto   &  \text{bar} \circ \Phi^b_c (x).
\end{eqnarray*}
This map is $\pi_1$-equivariant with respect to the homomorphism $f_*: \pi_1 (M) \longrightarrow \pi_1 (Z)$.  Thus the
straight line homotopy from $f: \U{M} \longrightarrow \U{Z}$ to $F_c: \U{M} \longrightarrow \U{Z}$ is similarly
equivariant, and descends to a homotopy between the maps $f, F_c: M \longrightarrow Z$.


\begin{lem} \label{Jac facts} For all $x \in M$, $ \text{Jac} \, F_c (x) \le (c/h)^n.$  Moreover,
\mbox{$(\text{Jac} \, F_c )  \longrightarrow 1$} in $L^1 (M, d \text{vol}_b)$ (and thus almost everywhere) as $c
\searrow h$. \end{lem} \begin{pf} Pick $x \in M$ and a positively oriented basis $\{ v_1, \ldots, v_n \} \subset
T_x M$.  Let $d\text{vol}_b$ (resp. $d\text{vol}_{g_0}$) be the volume form of $(M,b)$ (resp. $(Z,g_0)$).  Then
\begin{eqnarray*}
 \text{Jac} \, F_c  \cdot  d\text{vol}_b (v_1 , \ldots \, v_n )  &= &
         d\text{vol}_{g_0} (dF_c (v_1), \ldots, dF_c (v_n ))  \\
&=&  \Omega (d \Phi^b_c (v_1), \ldots , d \Phi^b_c (v_n)) \\ &\le &  \left( \frac{4n}{h^2} \right)^{n/2} \left(
\text{det} ( \langle d \Phi^b_c (v_i) \, , \,
            d \Phi^b_c (v_j ) \rangle_{L^2} )_{ij} \right)^{1/2} \quad \text{(Thm. \ref{calibration thm})}\\
&\le &  \left( \frac{4n}{h^2} \right)^{n/2} \, \left( \frac{c^2}{4n} \right)^{n/2}
            d \text{vol}_b (v_1, \ldots, v_n )  \quad \text{(Lemma \ref{volume inequality})}\\
           &=& \left( \frac{c}{h} \right)^n \ d\text{vol}_b (v_1, \ldots, v_n) .
\end{eqnarray*} By the inequalities at the end of Section \ref{BCG inequality section}, and the hypotheses of
Theorem \ref{equality prop}, it follows that 

\begin{eqnarray*}
\int_M \text{Jac} \, F_c \ d\text{vol}_b = \int_M (\Phi^b_c)^* \Omega & \le &
       \left( \frac{4n}{h^2} \right)^{n/2} \cdot \text{Vol}( \Phi^b_c)  \\
& \le &  \left( \frac{4n}{h^2} \right)^{n/2} \cdot \left( \frac{c^2}{4n} \right)^{n/2} \cdot \text{Vol}(M,b) \\
& = & \left( \frac{c}{h} \right)^n \cdot \text{deg}(f) \cdot \text{Vol}(Z,g_0) \\
&=& \left( \frac{c}{h} \right)^n \cdot \text{deg}(f) \cdot \left( \frac{4n}{h^2} \right)^{n/2} \cdot \text{Vol}(\Phi_0) \\
&=& \left( \frac{c}{h} \right)^n \cdot \text{deg}(f) \cdot \int_Z \Phi_0^* \Omega \\
&=& \left( \frac{c}{h} \right)^n \cdot \int_M (\Phi_0 \circ f)^* \Omega =
       \left( \frac{c}{h} \right)^n \cdot \int_M (\Phi^b_c)^* \Omega.
\end{eqnarray*}

From these inequalities we can conclude that 
$$\int_M \text{Jac} \, F_c \ d\text{vol}_b  \longrightarrow  \text{Vol}(M,b) = \int_M 1 \ d\text{vol}_b.$$
This limit and the fact that $\text{Jac} \, F_c \le (c/h)^n$ together imply that $(\text{Jac} \, F_c) \longrightarrow 1$ in $L^1 (M, d\text{vol}_b)$.
\end{pf}

We will now work with the lifted maps $F_c : \U{M} \longrightarrow \U{Z}$.  Pick a point $p \in \U{M}$.
For $\delta>0$ sufficiently small, let $B_p$ be a convex
$b$-metric ball of $b$-radius $\delta$ about $p$.  By reproducing the arguments of Lemmata 7.4 and 7.5 of
\cite{BCGlong}, one can show that the restricted maps $F_c |_{B_p}$ are uniformly $L$-Lipschitz for $c-h$
sufficiently small, where $L$ is independent of $p$.  (To apply the argument of \cite[Lem.7.5]{BCGlong} directly, it is necessary to restrict to a small ball.  $(\U{M},b)$ may not be uniformly Ahlfors regular, i.e. there may be a sequence of metric balls in $(\U{M},b)$ of fixed radius with volume going to zero.  Without uniform Ahlfors regularity the proof of \cite[Lem.7.5]{BCGlong} is not immediately valid for all of
$(\U{M},b)$.  Working locally is one way to circumvent this problem.)  

Recall that $(\U{Z}, g_0)$ has a basepoint $o \in \U{Z}$.
Then for each $c> h$ pick $g_c \in \text{Isom}^+ (\U{Z},g_0)$ such that $g_c \circ F_c |_{B_p} (p)=o$.  By the
Arzela-Ascoli theorem, there is a sequence $c_i \searrow h$ and an $L$-Lipschitz map $F^p : ( B_p , p)
\longrightarrow (\U{Z}, o)$ such that $g_c \circ F_{c_i} |_{B_p}$ converges uniformly to $F^p$.  By the arguments
of Lemmas 7.7 and 7.8 of \cite{BCGlong}, $F^p$ is in fact $1$-Lipschitz.

Pick a sequence of points $\{ p_i \} \subset \U{M}$ such that the neighborhoods $B_{p_i}$ form a locally finite cover of $\U{M}$.
Using a diagonalization argument, we may piece together the locally defined maps $F^{p_i}$ to prove

\begin{lem} \label{compact-open} There is a sequence $c_j \searrow h$ with the following property: for any
pre-compact open set $\mathcal{O} \subset \U{M}$, and any $\varepsilon >0$, there is a $J>0$ such that $j>J$
implies that the restricted maps $F_{c_j} : (\mathcal{O},b) \longrightarrow (\U{Z},g_0)$ are $L$-Lipschitz maps
satisfying the inequality $$d_{g_0}(F_{c_j} (x), F_{c_j}(y)) \le  d_b (x,y) + \varepsilon \quad \text{for any}
\quad x,y \in \mathcal{O}.$$
\end{lem}

In order to apply the Arzela-Ascoli theorem to the family of maps $\{ F_{c_j} \}_j$, we must first show that the
maps $F_{c_j}: (M,b) \longrightarrow (Z,g_0)$ do not leave every compact set of $Z$.

\begin{lem} \label{a-a prep} Pick a point $q \in M$.  The points $\{ F_{c_j} (q) \}_j \subset Z$ are contained in
a compact subset of $Z$. \end{lem} \begin{pf}
The map $f: M \longrightarrow Z$ is a proper map of finite degree.  A covering space argument shows that the image
group $f_* (\pi_1 (M)) \le \pi_1 (Z)$ must be of finite index.  The group $f_*(\pi_1(M))$ is therefore finitely
generated.  This implies the existence of a compact submanifold $C$ of $M$ containing $q$ such that the induced
map $$(f|_{C})_* : \pi_1 (C) \longrightarrow f_*(\pi_1(M))$$ is surjective.

Suppose that after passing to a subsequence the points $F_{c_j}(q)$ leave every compact subset of $M$.  Then for
large $j$, Lemma \ref{compact-open} implies that $F_{c_j}$ maps all of $C$ into a cusp of $(Z,g_0)$.  This implies
that the image of the homomorphism $$(F_{c_j}|_C )_* : \pi_1 (C) \longrightarrow \pi_1 (Z)$$ is a subgroup of
infinite index.  But $F_{c_j}$ is homotopic to $f$, implying that $$(F_{c_j}|_C )_* (\pi_1 (C)) = (f|_C)_* (\pi_1
(C)) = f_* (\pi_1(M)).$$ This is a contradiction. \end{pf}

We may now apply the Arzela-Ascoli theorem to conclude that, after passing to a subsequence, the maps $\{ F_{c_j}:
(M,b) \longrightarrow (Z,g_0) \}$ converge uniformly on compact sets to a $1$-Lipschitz map $F:(M,b) \longrightarrow
(Z, g_0)$ given locally by the maps $F^p : (B_p,p) \longrightarrow (\U{Z},o)$ described above.  The final and most important step of
the argument is to the following proposition.

\begin{prop} \label{jac is one} $\text{Jac} \, F = 1$ almost everywhere. \end{prop} 

\begin{pf} Since $F$ is $1$-Lipschitz, we know that $| \text{Jac} \, F | \le 1$ almost everywhere.  In search of a contradiction, let us suppose that $ \text{Jac} \, F  < 1$
on a set of strictly positive measure.  Then there is a $\gamma > 0$ and a compact submanifold $K \subset M$ such
that 
\begin{equation} \label{gap inequality}
\left( \int_K d\text{vol}_b \right) - \left( \int_K   \text{Jac} \, F \, d\text{vol}_b \right) > \gamma.
\end{equation}
(The letter $\gamma$ stands for ``gap.'')

Let us recall some notation from Section \ref{BCG inequality section}.  $M$ has a compact exhaustion $\{ M_i \}$, $L_i := \partial M_i$, and $\text{Vol}(L_i) \longrightarrow 0$.  Outside of a compact submanifold, $Z$ splits into the smooth product $Y \times [0,\infty)$, where each component of a slice $Y \times \{ t \}$ is the quotient of a horosphere in $\U{Z}$.  (Keep in mind that neither $L_i$ nor $Y$ is necessarily connected.)

Fix $i$ sufficiently large such that 
$$f(L_i) \subset Y \times [0, \infty) \subset Z, \quad K \subset M_i, \quad \text{and  } \text{Vol}(L_i ) < \frac{\gamma}{2}.$$
The following lemma is an exercise in degree theory whose proof is left to the reader.

\begin{lem} \label{degree lemma}
$f_* [ L_i] =  \text{deg} (f) \cdot [ Y ] \in H_{n-1} (Y \times [0,\infty)).$
\end{lem}

For each component $Y'$ of $Y$ there is a Riemannian covering space of $Z$ corresponding to $\pi_1 (Y') < \pi_1 (Z)$.  Take the disjoint union of these covering spaces, one for each component of $Y$, and call the resulting covering space $\hat{Z}$.  Let $\pi : \hat{Z} \longrightarrow Z$ denote the locally isometric covering map.  Note that $\hat{Z}$ is diffeomorphic to $Y \times \mathbb{R}$, and each component of $\hat{Z}$ is a quotient of the symmetric space $\U{Z}$ by a group of isometries with a global fixed point in the boundary at infinity $\partial \U{Z}$.

We chose $i$ so that $f(L_i) \subset Y \times [0,\infty) \subset Z$.  Therefore the restriction $f|_{L_i}$ lifts to a map $\widehat{ f|_{L_i}}: L_i \longrightarrow \hat{Z}.$  The map $F$ is homotopic to $f$, so in particular the restriction $F|_{L_i}$ also lifts to a map $\widehat{ F|_{L_i}} : L_i \longrightarrow \hat{Z}$.

Let $L_i \times (-\varepsilon, \varepsilon) \subset M$ be a tubular neighborhood of $L_i$, oriented so that $L_i \times \{ -\varepsilon/2 \} \subset M_i$.  Following the notation of Section \ref{cone section}, we define the cone of $\widehat{F|_{L_i}}$
$$\mathcal{C} : L_i \times [0, \varepsilon) \longrightarrow \hat{Z}$$
as follows.  For each point $p = \widehat{ F|_{L_i}}(x) \in \hat{Z}$ there is a unique geodesic ray beginning at $p$ which escapes down a finite volume end of $\hat{Z}$.  Define $\mathcal{C}(x,s)$ as the point at distance $\tan \left( \frac{s \pi}{2 \varepsilon} \right)$ along this ray.  Then by Corollary \ref{cone cor} we have the inequality
\begin{equation} \label{cone inequality 2}
\int_{L_i \times [0,\varepsilon)} | \text{Jac} \, \mathcal{C} | \, d \text{vol}_b \le
		\frac{1}{n-1} \, \int_{L_i} | \text{Jac} \, F|_{L_i}  | \, d \text{vol}_b \le \text{Vol}(L_i),
\end{equation}
where the final inequality comes from the fact that $F|_{L_i}$ is $1$-Lipschitz.

Fix $\tau > 0$ such that $\text{Vol}(Y \times [\tau, \infty) ) < \frac{\gamma}{2} \cdot \text{deg}(f)$.  Fix $\delta \in (0, \varepsilon)$ such that 
$$(\pi \circ \mathcal{C} ) (L_i \times \{ \delta \} ) \subset Y \times (\tau , \infty) \subset Z.$$
Define the compact submanifold $M^{\text{th}}_i := M_i \cup ( L_i \times [0,\delta]).$
(The letters th stand for ``thickening.'')  Define the compact manifold $Z_\tau$ as the complement of $Y \times (\tau, \infty) \subset Z$.  Define the piecewise smooth map $p: Z \longrightarrow Z$ as the identity on $Z_\tau$, and $p(y,t) := (y,\tau)$ for $(y,t) \in Y \times [\tau ,\infty)$. Notice that the Jacobian of $p$ is either $0$ or $1$ almost everywhere.  It is, in other words, a volume nonincreasing map (even though it is not globally Lipschitz).  See Figure \ref{cool picture} for a summary of this notation.

\begin{figure}
\psfrag{A}{$L_i = \partial M_i$ inside}
\psfrag{A'}{$\quad L_i \times (-\varepsilon, \varepsilon)$}
\psfrag{B}{$L_i \times \{\delta \} = \partial M^{\text{th}}_i$}
\psfrag{F}{$f,F$}
\psfrag{K}{$K$}
\psfrag{M}{$M$}
\psfrag{Z}{$Z$}
\psfrag{Y}{$Y \times \{ 0 \}$}
\psfrag{Y times tau}{$Y \times \{ \tau \}$}
\includegraphics{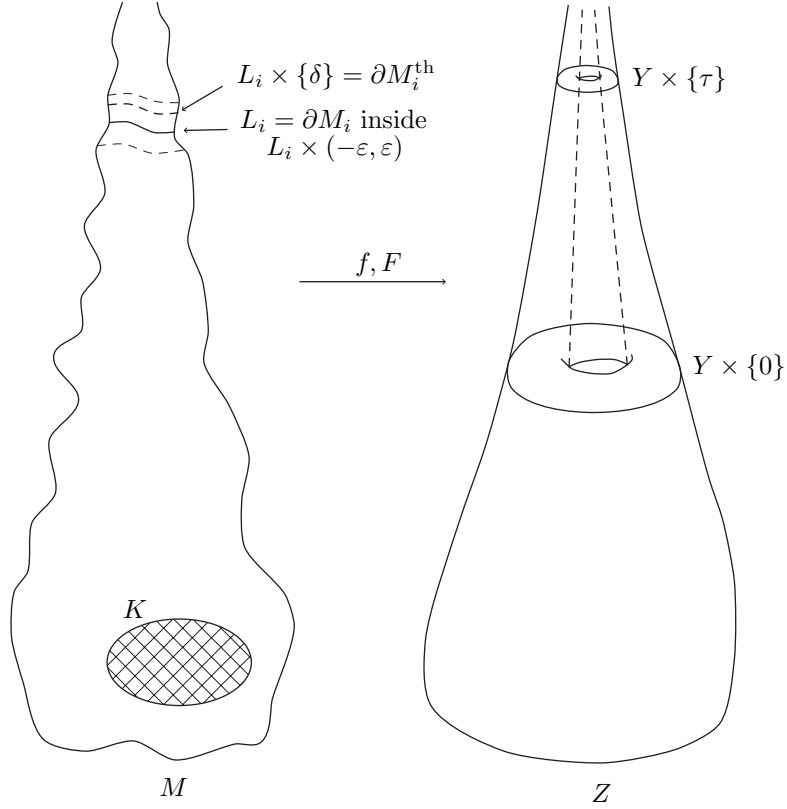}
\caption{A pictorial summary of the notation.} \label{cool picture}
\end{figure}

With this notation, we are prepared to define $G : M^{\text{th}}_i   \longrightarrow Z_\tau$ as follows: $G := p \circ F$ on $M_i$, and $G := p \circ \pi \circ \mathcal{C}$ on $L_i \times [0,\delta]$.  Note that $G$ is Lipschitz and boundary preserving.  Let $\text{deg}(G)$ denote the degree of $G$ as a boundary preserving map between compact manifolds.

\begin{lem} \label{deg G = deg f}
$\text{deg}(G) = \text{deg}(f)$.
\end{lem}

\emph{Proof of Lemma \ref{deg G = deg f}:}  By Lemma \ref{degree lemma},  $f_* [ L_i] =  \text{deg} (f) \cdot [ Y ] \in H_{n-1} (Y \times [0,\infty)).$  A standard fact of degree theory is that the degree of a boundary preserving map between compact manifolds equals the degree of the boundary map.  (This follows from the long exact homology sequence.)  In our notation, this implies the relation $G_* [L_i \times \{ \delta \}] = \text{deg}(G) \cdot [ \partial Z_\tau] \in H_{n-1} (\partial Z_\tau)$.

It is easy to see the map $G : M^{\text{th}}_i \longrightarrow Z$ is homotopic to the restriction of $f$ to $M^{\text{th}}_i$.  This implies the relation
$$G_* [ L_i \times \{ \delta \} ] = f_* [ L_i \times \{ \delta \} ] = \text{deg}(f) \cdot [Y]
	\in H_{n-1} (Y \times [0, \infty)),$$
where the restriction $G|_{L_i \times \{ \delta \}}$ is viewed as a map into $Y \times [0, \infty)$.  The proof of this lemma is completed by the commutative diagram

$$\xymatrix{
H_{n-1} (L_i \times \{ \delta \}) \ar[dr]_{G_*} \ar[rr]^{G_*} && H_{n-1}(Y \times [0, \infty)) \\
& H_{n-1}(\partial Z_\tau) \ar[ur]_{\text{incl}_*} &}$$
where the lower-right arrow is the isomorphism induced by inclusion.
\square

Consider the following estimate, with explanations below.
\begin{eqnarray*}
 \text{Vol}(Z_\tau)
&	= & \text{deg}(G) \, \int_{M^{\text{th}}_i} \text{Jac} \, G \ d\text{vol}_b \\
& = & \text{deg}(f) \int_{M_i} \text{Jac} \, ( p \circ F) \, d \text{vol}_b +
			\text{deg}(f) \int_{L_i \times [0, \delta]} \text{Jac} \, (p \circ \pi \circ \mathcal{C}) \, d \text{vol}_b \\
& < & \text{deg}(f) \cdot \left(  \text{Vol}(M_i) - \gamma \right)  + 
			\text{deg}(f) \int_{L_i \times [0,\delta]} | \text{Jac} \, \mathcal{C} | \, d \text{vol}_b \\
& < & \text{deg}(f) \cdot \text{Vol}(M) - \text{deg}(f) \cdot \gamma + \text{deg}(f) \cdot \text{Vol} (L_i)  \\
& < & \text{Vol}(Z) - \frac{\gamma}{2} \cdot \text{deg}(f).
\end{eqnarray*}
The first line is the degree formula for a boundary preserving Lipschitz map between compact smooth manifolds.  (Such a degree formula can be proved either by smooth approximation or by using \cite[Lem.4.1.25]{Fed}.)  The second line uses Lemma \ref{deg G = deg f} and the definition of $G$.  The third line uses the fact that $p: Z \longrightarrow Z$ is volume nonincreasing, inequality (\ref{gap inequality}), and the fact that $\pi : \hat{Z} \longrightarrow Z$ is a local isometry.  The fourth line uses inequality (\ref{cone inequality 2}).  The fifth line uses our choice of $i$ sufficiently large so that the volume of $L_i$ is less than $\frac{\gamma}{2}$.

Recall that we chose $\tau$ so that $\text{Vol}(Z) < \text{Vol}(Z_\tau) + \frac{\gamma}{2} \cdot \text{deg}(f)$.  This fact combined with the above string of inequalities yields a contradiction.  Therefore there does not exist a set of positive measure on which $\text{Jac} \, F < 1$.  This completes the proof of the proposition.

\end{pf}

We have now proven that $F: (M,b) \longrightarrow (Z,g_0)$ is a $1$-Lipschitz map such that $\text{Jac} \, F = 1$
almost everywhere.  We may now apply Theorem \ref{qr maps thm} to conclude that $F$ is a local homeomorphism and a
local isometry.  This implies that $(M,b)$ is a finite volume negatively curved locally symmetric space.  Since
$F$ is $1$-Lipschitz, it must take the cusps of $(M,b)$ out the cusps of $(Z,g_0)$.  Therefore $F$ is proper.  A
proper local homeomorphism is a covering map.  So we may conclude that $F:(M,b) \longrightarrow (Z,g_0)$ is a
locally isometric covering map.

It remains only to show that the straight line homotopy from $f$ to $F$ is proper.  For this it suffices to show
that an end of $(M,b)$ (which we now know to be a rank one cusp) is mapped to the same cusp of $(Z,g_0)$ under
both $f$ and $F$.  Suppose this is false.  Then there is an essential closed curve $\gamma$ in a cusp of $(M,b)$
which is mapped into different cusps of $(Z,g_0)$ under $f$ and $F$.  The homomorphism
  $$f_* = F_* : \pi_1(M) \longrightarrow \pi_1 (Z)$$
is injective.  Therefore the closed curves $F(\gamma)$ and $f(\gamma)$ are essential and freely homotopic in $Z$.
Since they lie in different cusps of $(Z,g_0)$, this violates the thick-thin decomposition of a negatively curved
locally symmetric space.  This yields a contradiction.  This concludes the proof of Theorem \ref{equality prop}.

\bibliography{unbdd.biblio} \end{document}